\documentclass[12pt,a4paper]{article}

\usepackage[french]{babel}
\usepackage[utf8]{inputenc}
\usepackage[T1]{fontenc}
\usepackage{amsmath}
\usepackage{array}
\usepackage{vmargin}
\usepackage{amssymb}
\usepackage{amsfonts}
\usepackage{enumerate}

\newtheorem{prop}{Proposition}[section]
\newtheorem{thm}[prop]{Theorem}
\newtheorem{cor}[prop]{Corollary}
\newtheorem{lemme}[prop]{Lemma}
\newtheorem{definition}[prop]{D\'efinition}
\newtheorem{remarque}[prop]{Remark}
\newtheorem{exemple}[prop]{Example}

\newtheorem{ex}{Exercise}
\newcommand{\qed}{\hfill \vrule height6pt width6ptdepth0pt}
\newenvironment{preuve}[1][]{\noindent {\it Proof #1}  : } {\hbox{~}\qed
\smallskip 

}
\newcommand{\field}[1]{\mathbb{#1}}
\newcommand{\R}{\field{R}}            
\newcommand{\N}{\field{N}}            
\newcommand{\C}{\field{C}}            

\newcommand{\SSS}{\field{S}}

\newcommand{\cB}{{\cal B}}

\newcommand{\cD}{{\cal D}}
\newcommand{\cE}{{\cal E}}

\newcommand{\cS}{{\cal S}}

\newcommand{\sgn}{{\rm sgn\,}}
\newcommand{\cotg}{{\rm cotg\,}}

\newcommand{\DATUM}{12-06-2015}              
\usepackage{fancyhdr}

\begin{document}
\title{ \textbf{Helffer-Sj\"ostrand formula for Unitary Operators.}}
\author{{\bf MBAREK Aiman}
}
\maketitle   
 ABSTRACT. The objective of this paper to give a formula for unitary operators that corresponds to the  Helffer-Sj\"ostrand formula for self-adjoint operators.

\setcounter{section}{0}
\section{Introduction}

First, we recall the usual formula of Helffer-Sj\"ostrand. If $f :\R\rightarrow\C$ a smooth function with compact support   and $A$ 
is a self-adjoint operator  on a Hilbert space $\mathcal{H}$, we have
$$f(A)=(2i\pi)^{-1}\int_{\C}\partial_{\bar{z}}f^{\C}(z)(z-A)^{-1}\,dz\wedge d{\bar{z}}.$$
Here  $f^\C$ is some almost analytic extension of $f$,  $dz\wedge d{\bar{z}}$ is the Lebesgue measure on the complex plane,  $\partial_{\bar{z}}=\frac 12(\partial_x+i\partial_y)$ for $z=x+iy,$  and $f(A)$  is given by the functional calculus for the self-adjoint operator $A$.
 This formula also holds true for a larger class of functions $f$ having some prescribed behaviour at infinity (See.  \{\cite{DG},\cite{GJ},...\}). For instance (See \cite{DG}), one can require, for $\rho<0$, that  $f\in\mathcal{C}^\infty(\R)$ such that
$$\forall k\in\N, \ \sup_{t\in\R}\langle t\rangle^{-\rho+k}|f^{(k)}(t)|<+\infty.$$

The  Helffer-Sj\"ostrand formula is  extensively used in many different works, for example (\cite{DG},\cite{GJ},\cite{BG},\cite{Ca},\cite{GN},...). A first application is an expansion of commutators of the following type.
Given $f$ a smooth function  on $\R$ with compact support, $A$  a self-adjoint operator  and $B$  a bounded operator satisfying certain properties  on a Hilbert space $\mathcal{H}$, one has
$$[f(A),B]=f'(A)[A,B]+R,$$
where  $[A,B]$ is the commutator of $A$ and $B$, $f(A)$ is given by the functional calculus for $A$ and $R$  is a rest that  has  better properties in some sense. One can generalize this formula with  iterated  commutators $([A,[A,B]],[A,[A,[A,B]]],...)$ to get a Taylor-type formula.\\

A second application is provided in \cite{GN}. Given a self-adjoint operator $A$ on $\mathcal{H}$ and  a smooth function $f$  on $\R$ with compact support, one gets a  control on the norm $\|f(A)\|_{\mathcal{B}(\mathcal{H})}$  in terms of the norm of the resolvent of $A$,                 
          $\|(A-z_0)^{-1}\|_{\mathcal{B}(\mathcal{H})}$, for some $z_0\in \C\setminus \R$. There are also results of this kind where $\mathcal{B}(\mathcal{H})$  is replaced by a Schatten class norm.\\

Let us mention a third application in the context of linear PDE. For $\zeta>0$, let $V$ be the multiplication operator by a function $V:\R^d\rightarrow\R,x\mapsto V(x)$ such that $x\mapsto|\langle x\rangle^{\zeta}.V|$ is 
 bounded, $\langle Q\rangle$ be the multiplication operator by the map $x\mapsto\langle x\rangle=(1+|x|^2)^{\frac 12}$, $f$ be a smooth function  on $\R$ with compact support  and $-\Delta$ be the positive Laplacian on $\R^{d}$. Let $H_1 =H_0+V=-\Delta  + V$, the self-adjoint realization in $\mathcal{L}^2(\R^d)$.  By the Helffer-Sj\"ostrand formula one can show that\\ 
$$\big((f(H_1)-f(H_0))\langle Q\rangle^{\epsilon}\langle H_0\rangle^{\alpha})_{0\leq\epsilon<\zeta,0\leq\alpha<1}$$ is a family of  compact operators. \\

Our objective is to find a kind of  Helffer-Sj\"ostrand  formula for unitary operators $U$ and for smooth functions defined on the 1-dimensional sphere $\SSS^1$. Since $U$ is a unitary operator, its spectrum is contained in $\SSS^1$ and a functional calculus  is well defined. Here we shall take a smooth function $f:\SSS^1\rightarrow\C$  with compact support  in the interior of $\SSS^1\setminus\{1\}$. The latter condition corresponds to the condition of compact support for the Helffer-Sj\"ostrand formula. For a unitary operator $U$ on some  Hilbert space $\mathcal{H},$ we shall derive the formula
\begin{equation}\label{uniH-J}
f(U)=(2i\pi)^{-1}\int_{\C}\partial_{\bar{z}}f^{\C}_{\SSS^1}(z)(z-U)^{-1}\,dz\wedge d{\bar{z}},
\end{equation}
where  $f^\C_{\SSS^1}$ is some almost analytic extension of $f$, $dz\wedge d{\bar{z}}$ is the Lebesgue measure on the complex plane and $f(U)$  is given by the functional calculus for $U$. We  will show the formula (\ref{uniH-J})  in Theorem \ref{ssss}. We expect that our formula (\ref{uniH-J}) holds true under a weaker assumption on the behaviour of $f$ near 1.\\

The tools used in this paper come from complex analysis (Cayley transform...), the theory of self-adjoint, unitary and normal operators and the functional
 calculus  for the  corresponding operators.\\

The paper is organized as follows. In Section 2, we formulate the theorem of Helffer-Sj\"ostrand  and  we give a  slightly different proof from that of \cite{hs}. In Section 3, we recall some properties on the Cayley transform and we  prove some complex analysis results that will be used in the proof of our main theorem. In Section 4, we prove the main theorem of this paper. The paper ends with a paragraph of notation.\\

 I would like  to thank  S. Golénia, who drew our attention to the fact that a  formula of the type (\ref{uniH-J}) should exist and would be very interesting and useful.  I also thank T. Jecko  for guiding, abetting, counseling me throughout the preparation of this paper.

\newpage
\section{Usual Helffer-Sj\"ostrand formula}
We state the theorem of Helffer Sj\"ostrand  and  we will give another proof other than in the article \cite{hs}. Our proof is based
 on results in \cite{DG} and results of complex analysis.

First we will start with constructing almost analytic extensions.

Let $k : \C \longrightarrow \C$ be a map,  we denote by $supp k$  the support of $k$.
\begin{prop}\cite{DG}\label{complex} 
Let $f$ a smooth function  on $\R$ with compact support. Then there exists a smooth function $f^\C:\C\rightarrow\C$, called   
 an almost analytic extension of $f$, such that  
 there exists $0<C<1$ such that for all  $l\in \N$ there exists $C_l\geq0$ such that, \\
\begin{align}
   f^{\C}|_{\R}=f, supp f^{\C}\subset\{x+i y;x\in supp f , |y|\leq C\},\label{cplx}\\ 
| \partial_{\overline z}f^{\C}(z)|\leq C_l|Im(z)|^l. \label{cpl}
\end{align}
where $Im(z)$ is the imaginary part of $z$.
\end{prop}

\begin{preuve}
we follow the argument of \cite{DG}, checking that we can ensure the formula (\ref{cplx}) with $C <1.$

 Let $M_n=\sup_{0\leq k\leq n}\sup_{x\in\R}|f^{(k)}(x)|$ for all $n\in\N$. 
 Let $\chi\in\mathcal{C}^{\infty}(\R)$ such that $\chi=1$ on $[-1/2,1/2]$ and $\chi=0$ on $\R\setminus{[-1,1]}$. We choose
\begin{equation}\label{sume}
f^{\C}(z)=f^{\C}(x+iy)=\displaystyle\sum_{n=0}^{+\infty}\frac{(iy)^n}{n!}f^{(n)}(x)\chi\big(\frac{y}{T_n}\big),
\end{equation}
for a  decreasing and positive sequence $(T_n)_n$ satisfying, $T_0<1$ and 
\begin{equation}\label{Tn}
\forall n\in\N,\ \ \ \  T_n M_{2n}\leq 2^{-n}.
\end{equation}
For example  we can take $(T_n)_n$ defined by  $0<T_0<1$ and $T_n=min(T_{n-1},\frac{1}{2^n M_{2n}}),$ for $ n\geq1$.
 With this choice we will see that our sum (\ref{sume}) is uniformly convergent. Then, for $x$ real,
$$f^{\C}(x)=\displaystyle\sum_{n=0}^{+\infty}\frac{(i0)^n}{n!}f^{(n)}(x)\chi\big(\frac{0}{T_n}\big)=\frac{(i0)^0}{0!}f(x)\chi\big(\frac{0}{T_0}\big)=f(x).$$

We deduce that $f^{\C}$ is an extension of $f$ to $\C.$\\
Let $C=T_0$. Let $H_{f,C}:=\{x+i y;x\in suppf , |y|\leq C\}$. If $x_0+iy_0\in\C$ such that the distance  $d(x_0+iy_0,H_{f,C})>0$ then $x_0\notin supp f$ 
or $|y_0|>C$. If $x_0\notin supp f$ therefore $f^{(n)}(x_0)=0, n\geq0$  for all $n\in\N,$ and $f^{\C}(x_0+iy_0)=0$. If $|y_0|>C$ then, for all $n\in\N, |y_0|>C\geq T_n$, since  $(T_n)_n$ decreases. Thus, for all $n, \chi(\frac{y_0}{T_n})=0$ and 
 $f^{\C}(x_0+iy_0)=0$. We get (\ref{cplx}) \\

To show that our sum (\ref{sume}) exists and belongs to $\mathcal{C}^{\infty}(\C)$, we simply have to show that, for  $(p,q)\in\N^2$,

\begin{equation}\label{sumd}
\displaystyle\sum_{n\in\N}^{}\partial_x^p\partial_y^q\bigg(\frac{(iy)^n}{n!}f^{(n)}(x)\chi\big(\frac{y}{T_n}\big)\bigg)
\end{equation}

  is uniformly convergent.\\

Let $(p,q)\in\N^2$, take  $n\geq max(p,q)$. Using the Leibnitz formula  we have

\begin{equation}
\begin{array}{cl}
\bigg|\partial_x^p\partial_y^q\bigg(\frac{(iy)^n}{n!}f^{(n)}(x)\chi\big(\frac{y}{T_n}\big)\bigg)\bigg|&\\

=\bigg|i^n(\partial_x^{p+n}f)(x)\displaystyle\sum_{r=0}^{q}\frac{C^r_q}{(n-q+r)! \ T_n^r}y^{n-q+r}\chi^{(r)}\big(\frac{y}{T_n}\big)\bigg|.&
\end{array}
\end{equation}

where $C^r_q=\frac{q!}{r!(q-r)!}, 0\leq r\leq q$. For all $0\leq r\leq q\leq n, \chi^{(r)}(\frac{\cdot}{T_n})$ is supported in $\{ 0\leq|\cdot|\leq T_n\}$.\\
Let $x,y\in\R$, we have

\begin{equation}
\begin{array}{ll}
\bigg|\partial_x^p\partial_y^q\bigg(\frac{(iy)^n}{n!}f^{(n)}(x)\chi\big(\frac{y}{T_n}\big)\bigg)\bigg|

&\leq M_{p+n}\sum_{r=0}^{q}\frac{C^r_q}{(n-q+r)!T_n^r}{|y|}^{n-q+r}\big|\chi^{(r)}\big(\frac{y}{T_n}\big)\big|

\\&\leq M_{p+n}T_n^{n-q}\displaystyle\sum_{r=0}^{q}\frac{C^r_q}{(n-q+r)!}\|\chi^{(r)}\|_{\infty}

\\&\leq D_q M_{p+n}T_n^{n-q}.

\end{array}
\end{equation}

As $T_n<1$, we have $T_n^{n-q}M_{p+n}\leq T_n M_{2n}\leq 2^{-n},$ by (\ref{Tn}).

Thus 

$$\sup_{x,y}\bigg|\partial_x^p\partial_y^q\bigg(\frac{(iy)^n}{n!}f^{(n)}(x)\chi\big(\frac{y}{T_n}\big)\bigg)\bigg|\leq D_q 2^{-n}.$$

This yields the normal convergence of the series (\ref{sumd}). In particular, $f^{\C}\in\mathcal{C}^{\infty}(\C)$.\\

It remains to show (\ref{cpl}). For $x,y\in\R$
\begin{equation}
\begin{array}{ll}
2\partial_{\bar{z}}f^{\C}(z)

&=(\partial_x+i\partial_y)f^{\C}(x+iy)\\

&=\displaystyle\sum_{n=0}^{+\infty}\bigg(\frac{(iy)^n}{n!} f^{(n+1)}(x)\chi\big(\frac{y}{T_n}\big)+i\frac{i^n}{(n-1)!}f^{(n)}(x)y^{n-1}\chi\big(\frac{y}{T_n}\big)\\
&\ \ \ \ \ \ \ \ \ \ \ \ \ \ \ \ \ \ \ \ \ \ \ \ \ \ \ +i\frac{(iy)^n}{n!}f^{(n)}(x)\frac{1}{T_n}\chi'\big(\frac{y}{T_n}\big)\bigg)\\

&=f'(x)\chi\big(\frac{y}{T_0}\big)+\displaystyle\sum_{n=1}^{+\infty}\frac{(iy)^n}{n!} f^{(n+1)}(x)\chi\big(\frac{y}{T_n}\big)
+i\displaystyle\sum_{n=0}^{+\infty}\frac{f^{(n)}(x)}{(n-1)!}\bigg(i^n ny^{n-1}\chi\big(\frac{y}{T_n}\big)\\

&\ \ \ \ \ \ \ \ \ \ \ \ \ \ \ \ \ \ \ \ \ \ \ \ \ +(iy)^n\frac{1}{T_n}\chi'\big(\frac{y}{T_n}\big)\bigg)\\

&=f'(x)\chi\big(\frac{y}{T_0}\big)+\displaystyle\sum_{n=1}^{+\infty}\frac{(iy)^n}{n!} f^{(n+1)}(x)\bigg(\chi\big(\frac{y}{T_n}\big)-\chi\big(\frac{y}{T_{n+1}}\big)\bigg)\\

&\ \ \ \ \ \ \ \ \ \ \ \ \ \ \ \ \ \ \ \ \ \ \ \ \ +\displaystyle\sum_{n=0}^{+\infty}\frac{f^{(n)}(x)}{n! T_n}i^{n+1}y^n\chi'\big(\frac{y}{T_n}\big).
\end{array}
\end{equation}

Let $l\in\N$. Denoting by $\mathbf{1}_{I}$  the  characteristic function of $I$, we have, for $y\neq 0$,
\begin{equation}
\begin{array}{ll}
2\big| |y|^{-l}\partial_{\bar{z}}\widetilde{f}(z)\big|

&\leq (\frac{T_0}{2})^{-l}M_1\|\chi\|_\infty+\displaystyle\sum_{n=1}^{+\infty}M_{n+1}|y|^{n-l}2\|\chi\|_{\infty}\mathbf{1}_{\{T_{n+1}/2\leq|.|\leq T_n\}}(y)\\

&+\displaystyle\sum_{n=0}^{+\infty}\frac{M_n}{T_n}|y|^{n-l}\|\chi'\|_{\infty}\mathbf{1}_{\{T_n/2\leq|.|\leq T_n\}}(y)\\

&\leq C+C'\displaystyle\sum_{n=1}^{l}M_{n+1}\bigg(\frac{T_{n+1}}{2}\bigg)^{n-l}+C''\displaystyle\sum_{n=l+1}^{+\infty}M_{n+1}(T_n)^{n-l}

\\&+C'''\displaystyle\sum_{n=0}^{l}\frac{M_n}{T_n}\bigg(\frac{T_n}{2}\bigg)^{n-l}+C''''\displaystyle\sum_{n=l+1}^{+\infty}\frac{M_n}{T_n}(T_n)^{n-l}

\\&\leq C_l'+D\displaystyle\sum_{n=l+2}^{+\infty}\big(M_n+M_{n+1}\big)T_n.
\end{array}
\end{equation}

By (\ref{Tn}) and the definition of $(M_n)_n$  we have $\big(M_n+M_{n+1}\big)T_n\leq 2 T_n M_{2n} \leq 2.2^{-n}$. This yields (\ref{cpl}).
\end{preuve}

\begin{prop}\cite{H}\label{H\"orm}

Let $\omega$ be an open set in the complex plane $\C$ and $u\in\mathcal{C}^1(\omega)$. For all $ \xi\in\omega,$  the second integral in (\ref{integra}) below  is convergent. Moreover, for all $ \xi\in\omega$, 

\begin{equation}\label{integra}
u(\xi)=(2i\pi)^{-1}\bigg\{\int_{\partial\omega}\frac{u(z)}{z-\xi}\,dz+\int_{\omega}\frac{\partial_{\bar{z}}u(z)}{z-\xi}\,dz\wedge d{\bar{z}}\bigg\}.
\end{equation}
\end{prop}

\begin{cor}\label{compac}

If $u\in\mathcal{C}^{\infty}(\C)$ with compact support  in $\omega$  then, for all $\xi\in\omega,$

\begin{equation}\label{compactif}
u(\xi)=(2i\pi)^{-1}\int_{\omega}\frac{\partial_{\bar{z}}u(z)}{z-\xi}\,dz\wedge d{\bar{z}}.
\end{equation}

\end{cor}

\begin{preuve}
Just apply Proposition (\ref{H\"orm}) and use the fact that $u$ is zero on $\partial\omega.$
\end{preuve}

We will recall some results on the theory of normal operators. 
Let  $ f\in\mathcal{C}^{\infty}(\R)$ with compact support  $K$ and  let $N$ is a normal operator  on a  Hilbert space $\mathcal{H}$.
Denote by $\mathcal{B}(\mathcal{H})$  the set of bounded operators on $\mathcal{H}$ and $f(N)$ in $\mathcal{B}(\mathcal{H})$ given by the functional calculus for $N.$

\begin{prop}\cite{Co}\label{resolvo}
For  $z$ outside $\sigma(N)$, the spectrum of the operator $N$, we have

$$\|(z-N)^{-1}\|_{\mathcal{B}(\mathcal{H})}=\frac{1}{d\big(z,\sigma(N)\big)}.$$

\end{prop}

\begin{prop}\cite{Co}\label{conv}
Let  $K$ a compact of $\C$. Let, for $n\in\N, f_n\in\mathcal{C}^{\infty}(\C)$    with $supp f_n\subset K$ and $ f\in\mathcal{C}^{\infty}(\C)$ with $supp f\subset K$ such that $(f_n)_n$ converges  to $f$, uniformly on $K$.
  Then  $(f_n(N))_n$ converges to $f(N)$  for the operator norm.
\end{prop}

\begin{preuve}

By combining the lemma (1.9) page $257$ and Theorem (4.7) page $321$ in the book \cite{Co}, we found the following inequality

$$\|f_n(N)-f(N)\|_{\mathcal{B}(\mathcal{H})}\leq\|f_n-f\|_{\infty,K}:=\sup_{x\in K}|f_n(x)-f(x)|.$$

Therefore by passing to the limit we conclude

$$\lim_{n\rightarrow+\infty}\|f_n(N)-f(N)\|_{\mathcal{B}(\mathcal{H})}=0.$$
\end{preuve}
If we replace $u$ by an almost analytic extension $f^\C$  satisfying (\ref{compac}) and take $\omega=\C$, the generalized integral in (\ref{compactif}) converges uniformly, w.r.t. $\xi$, as proved in the following proposition.

For a map $k : \C \longrightarrow \C$  and $D\subset\C$, we define

\begin{equation}\label{D}
\|k\|_{\infty,D}:=\sup_{x\in D}|k(x)|.
\end{equation}

\begin{prop}\label{convergence}
Let $f\in\mathcal{C}^{\infty}(\R)$ be a function with compact support  $K$ and $f^{\C}$  be an almost analytic extension of $f$ given by Proposition \ref{complex}. 
We have the following  uniform convergence on $K$

\begin{equation}\label{norminf}
\lim_{n\rightarrow+\infty}\bigg\|\int_{|Im(z)|>1/n}\partial_{\bar{z}}f^{\C}(z)(z-\cdot)^{-1}\,dz\wedge d{\bar{z}}-f(\cdot)\bigg\|_{\infty,K}=0.
\end{equation}

\end{prop}

\begin{preuve}
Let $f\in\mathcal{C}^{\infty}(\R)$ be a function with compact support  K. By Proposition \ref{complex}, $f^{\C}$  
is an almost analytic extension of $f$ such that $supp f^{\C}\subset K^{\C}$, where \\$K^{\C}:=\{(x,y)\in \R^2 ; x\in K ,|y|\leq C< 1\}$.
 Let $\omega$ an open in the complex plane $\C$ that contains $K^{\C}$.\\ From the corollary \ref{compac} with $u=f^{\C}$, we have, for all $\xi\in K$,

\begin{equation}\label{integ}
f(\xi)=(2i\pi)^{-1}\int_{\C}\frac{\partial_{\bar{z}}f^{\C}(z)}{z-\xi}\,dz\wedge d\bar{z} .
\end{equation}

Using the  (\ref{compactif}) with $l = 1$, we have, for $\xi\in K$,

\begin{equation}\label{etoile}
\int_{K^{\C}}|\partial_{\bar{z}}f^{\C}(z)(z-\xi)^{-1}|\,dz\wedge d{\bar{z}}\leq C_1\int_{K^{\C}}|Im(z)||Im(z)|^{-1}\,dz\wedge d{\bar{z}}<+\infty .
\end{equation}

Thus this integral (\ref{integ}) is absolutely convergent. Now we will show the uniform convergence on $K$, let $\xi\in K$

\begin{equation}
\begin{array}{ll}

\bigg|\int_{|Im(z)|>1/n}\partial_{\bar{z}}f^{\C}(z)(z-\xi)^{-1}\,dz\wedge d{\bar{z}}-f(\xi)\bigg|&\\

=\bigg|\int_{|Im(z)|<1/n}\partial_{\bar{z}}f^{\C}(z)(z-\xi)^{-1}\,dz\wedge d{\bar{z}}\bigg|&\\

\leq C_1\int_{|Im(z)|<1/n}|Im(z)||Im(z)|^{-1}\,dz\wedge d{\bar{z}}&\\

\leq \frac{C}{n} \ \ \mbox{\ \ \ , where C is independent of }\xi.&

\end{array}
\end{equation}

Therefore

$$\sup_{\xi\in K}\bigg|\int_{|Im(z)|>1/n}\partial_{\bar{z}}f^{\C}(z)(z-\xi)^{-1}\,dz\wedge d{\bar{z}}-f(\xi)\bigg|< \frac{C}{n},$$

proving (\ref{norminf}).
\end{preuve}

We are able to reprove  Helffer-Sj\"ostrand theorem.
\begin{thm}\cite{hs}
Let $f\in\mathcal{C}^{\infty}(\R)$ be a function with compact support  $K$ and $f^{\C}$  be an almost analytic extension of $f$ given by Proposition \ref{complex}.
 Let $A$ be a self-adjoint operator  on a  Hilbert space $\mathcal{H}$. The integral\\
\begin{equation}\label{doubleetoile}
\int_{\C}\|\partial_{\bar{z}}f^{\C}(z)(z-A)^{-1}\|_{\mathcal{B}(\mathcal{H})}\,dz\wedge d{\bar{z}}
\end{equation}

converges, the integral in the following formula converges in operator norm in $\mathcal{B}(\mathcal{H})$, and we have

\begin{equation}\label{Sjos}
f(A)=(2i\pi)^{-1}\int_{\C}\partial_{\bar{z}}f^{\C}(z)(z-A)^{-1}\,dz\wedge d{\bar{z}}.
\end{equation}

\end{thm}

\begin{preuve}
By Proposition \ref{resolvo}, (\ref{etoile}) holds true with $\xi$ replaced by $A$. This proves the convergence of  (\ref{doubleetoile}).   In particular, the integral in (\ref{Sjos}) converges in the operator norm of $\mathcal{B}(\mathcal{H})$. By(\ref{norminf}) and  Proposition \ref{conv} we get \eqref{Sjos} .
\end{preuve}

\section{Cayley transform}
In this section, we give some properties on the Cayley transform and we prove a known result in complex analysis that will be used in the proof of our main theorem.

 We will need some estimates on the Cayley transform on specific regions.

\begin{definition}
The Cayley transform is the map
\begin{equation}
\begin{array}{clcl}
\psi&:\C\setminus\{i\}&\longrightarrow&\C\setminus\{1\}\\&\ \ \ \ \ z&\longmapsto&\psi(z)=\frac{z+i}{z-i} .
\end{array}
\end{equation}
\end{definition}
We denote by  $\SSS^1$ is the unit sphere in $\C$, $D(0,1)$ the open disk with center the origin  and of radius $1$, $ \overline{D(0,1)}$ the closure of  $D(0,1)$, $\C^+=\{z\in\C ; Im(z)>0\}$ and  $\C^-=\{z\in\C ; Im(z)<0\},$ where $Im(z)$ is the imaginary part of $z$. 
\begin{prop}\cite{H}
$\psi$ is an analytic, bijective function and $\psi^{-1}$ is given by $\psi^{-1}(\xi)=i\frac{\xi+1}{\xi-1},$ for all $\xi\in\C\setminus\{1\}.$ Furthermore, 
$\psi(\R)=\SSS^1\setminus\{1\}, \psi(\C^-)=D(0,1), \psi(\C^+\setminus\{i\})\\=\C\setminus\overline{D(0,1)}$,
$\psi(0)= -1 ,\psi(-1)=i , \lim_{z\rightarrow\infty}\psi(z)=1, \lim_{z\rightarrow i}|\psi(z)|=\infty .$\\
\end{prop}

Let  $a,b,c\in\R$ such that $a<b$ and $0\leq c<1$. We define
$$\Omega:=\psi([a,b]\times[-c,c]).$$

\begin{lemme}\label{majpsi}
There exists $C>0,$ such that
\begin{equation}\label{psinfn}
\|\partial_x{(\psi^{-1})}\|_{\infty,\Omega}\leq C, \|\partial_y{(\psi^{-1})}\|_{\infty,\Omega}\leq C,
\end{equation}

and,  for all $\xi\in\Omega,$ 
$$|Im(\psi^{-1}(\xi))|\leq C d(\xi, \SSS^1).$$
where $\|.\|_{\infty,\Omega}$ defined in (\ref{D}) and $d(\xi, \SSS^1)$  is the distance between $\xi$ and $\SSS^1$.

\end{lemme}

\begin{preuve}
Let $\xi\in\Omega$ . Then $|Im\big(\psi^{-1}(\xi)\big)|=\bigg|\frac{|\xi|^2-1}{|\xi-1|^2}\bigg|$,
$$\frac{|Im\big(\psi^{-1}(\xi)\big)|}{d(\xi, \SSS^1)}=\frac{|\xi|+1}{|\xi-1|^2}.$$
Call $C_\Omega=d(1,\Omega)>0, d_\Omega=\sup_{\xi\in\Omega}|\xi|>0$. We have, for all $\xi\in\Omega$,

$$C_1:=\frac{1}{1+d_\Omega}\leq\frac{1}{|\xi|+1}\leq\frac{|Im\big(\psi^{-1}(\xi)\big)|}{d(\xi, \SSS^1)}=
\frac{|\xi|+1}{|\xi-1|^2}\leq\frac{d_\Omega+1}{C_\Omega^2}=:C_2.$$

$\partial_x{(\psi^{-1})}=\frac{-2i}{(\xi-1)^2}$ and $\partial_y{(\psi^{-1})}=\frac{2}{(\xi-1)^2}$. 
 Since, for all  $\xi\in\Omega, |\xi-1|\geq C_{\Omega}>0,$ we get (\ref{psinfn}).
\end{preuve}

\begin{lemme}\label{derive}
Let $h$ be an analytic function on $\C$ and $g$ be a smooth function on $\C\simeq\R^2$. We have 
$$\forall z\in\C, \ \ \  \partial_{\bar{z}} (g\circ h)(z)=(\partial_{\bar{z}}g)(h(z))\partial_{\bar{z}}\overline{h(z)} .$$

\end{lemme}
\begin{preuve}
Writing $z=x+iy$ and  $h(z)= h_1(x,y)+i h_2(x,y)$, we have
\begin{equation}
\begin{array}{cl}
2\partial_{\bar{z}}(g\circ h)(z)&=\partial_x(g\circ h)(z)+i \partial_y(g\circ h)(z)\\&=(\partial_x g)(h(z))\partial_x h_1(x,y)+ (\partial_y g)( h(z))\partial_x h_2(x,y)\\

&+i\big((\partial_x g)(h(z))\partial_y h_1(x,y)+ (\partial_y g)( h(z))\partial_y h_2(x,y)\big)\\

&=((\partial_z+\partial_{\bar{z}}) g)(h(z))\partial_x h_1(x,y)+ (\frac{\partial_{\bar{z}}-\partial_z}{i}) g)( h(z))\partial_x h_2(x,y)\\

&+i\big(((\partial_z+\partial_{\bar{z}}) g)(h(z))\partial_y h_1(x,y)+ ((\frac{\partial_{\bar{z}}-\partial_z}{i}) g)( h(z))\partial_y h_2(x,y))\big)\\

&=(\partial_z g)(h(z))\big(\partial_x h_1(x,y)+i\partial_x h_2(x,y)+i\partial_y h_1(x,y)-\partial_y h_2(x,y)\big)\\
&+(\partial_{\bar{z}}g)(h(z))\big(\partial_x h_1(x,y)-i\partial_x h_2(x,y)+ i\partial_y h_1(x,y)+\partial_y h_2(x,y)\big).
\end{array}
\end{equation}

 As $h$ is analytic, $h$ satisfies the Cauchy Riemann relations $\partial_x h_1(x,y)=\partial_y h_2(x,y)$ and $\partial_y h_1(x,y)=-\partial_x h_2(x,y)$, therefore
\begin{equation}
\begin{array}{cl}
2\partial_{\bar{z}}(g\circ h)(z)

&=2(\partial_{\bar{z}}g)(h(z))\big(\partial_x h_1(x,y)+i\partial_y h_1(x,y)\big)\\

&=2(\partial_{\bar{z}}g)(h(z))\big(\frac{\partial_x+i\partial_y}{2} h_1(x,y)+\frac{\partial_x+i\partial_y}{2}h_1(x,y)\big )\\

&=2(\partial_{\bar{z}}g)(h(z))\big(\frac{\partial_x+i\partial_y}{2} h_1(x,y)+\frac{\partial_y -i\partial_x}{2} h_2(x,y)\big)\\

&=2(\partial_{\bar{z}}g)(h(z))\big(\frac{\partial_x+i\partial_y}{2} h_1(x,y)-i\frac{\partial_x+i\partial_y}{2}h_2(x,y)\big)\\

&=2(\partial_{\bar{z}}g)(h(z))\partial_{\bar{z}} \overline{h(z)} .
\end{array}
\end{equation}
\end{preuve}

\section{Helffer-Sj\"ostrand formula for Unitary Operators}
We want a formula similar to (\ref{Sjos}),  replacing  the self-adjoint operator $A$ by a unitary operator $U$.

We start with some reminders on unitary operators.
\begin{definition}
$U$ is called a unitary operator on a  Hilbert space $\mathcal{H}$,  if $$U U^* =U^*U=I.$$
\end{definition}

\begin{prop}
The spectrum of a unitary operator $U$  on a Hilbert space $\mathcal{H}$  is included in $\SSS^1$ .
\end{prop}
\begin{preuve}
Let  $z\in\C$ such that $|z|>1$. We have $\|z^{-1} U\|<1$, then $(I-z^{-1} U)$ is invertible and its inverse is bounded and it is given by the following series $\displaystyle\sum_{n=0}^{+\infty}(z^{-1}U)^n$. Therefore $(z-U)z^{-1}(I-z^{-1}U)^{-1}=\sum_{n=0}^{\infty}(z^{-1}U)^n-\sum_{n=0}^{\infty}Uz^{-1}(z^{-1}U)^n=I=z^{-1}(I-z^{-1}U)^{-1}(z-U),$
$(z-U)$ is invertible, and  $(z-U)^{-1}=z^{-1}(I-z^{-1}U)^{-1}.$ 

 Let now $z\in\C$ such that $|z|<1$, we have $|z U^*|<1,$ then

$$(z-U)U^*(z U^*-I)^{-1}=\sum_{n=0}^{\infty}zU^*(zU^*)^n-\sum_{n=0}^{\infty}(zU^*)^n=I=(z-U)U^*(z U^*-I)^{-1},$$
$(z-U)$ is invertible, and  $(z-U)^{-1}=U^*(z U^*-I)^{-1}.$\\
 This shows that the spectrum of $U$  is included in $\SSS^1$.
\end{preuve}

~\\
We are able to demonstrate the main theorem of this paper.

\begin{thm}\label{ssss}
 Let $f\in\mathcal{C}^\infty(\SSS^1)$ supported away from 1 and $U$ be a unitary operator on a  Hilbert space $\mathcal{H}$. Let $(f\circ\psi)^{\C}$ 
 be an almost analytic extension of $f\circ\psi$ given by (\ref{complex}) such that  $(f\circ\psi)^{\C}$ is supported in $\{x+i y;x\in supp (f\circ\psi),
 |y|\leq C\}, 0<C<1$. Let $f^{\C}_{\SSS^1}=(f\circ\psi)^{\C}\circ\psi^{-1}$. The integral
\begin{equation}\label{int}
\int_{\C}\|\partial_{\bar{z}}f^{\C}_{\SSS^1}(z)(z-U)^{-1}\|_{\mathcal{B}(\mathcal{H})}\,dz\wedge d{\bar{z}}
\end{equation}
 converges, the integral in the following formula converges in operator norm  in $\mathcal{B}(\mathcal{H})$, and we have
\begin{equation}\label{uniHJ}
f(U)=(2i\pi)^{-1}\int_{\C}\partial_{\bar{z}}f^{\C}_{\SSS^1}(z)(z-U)^{-1}\,dz\wedge d{\bar{z}}.
\end{equation}
\end{thm}

\begin{preuve}
Let $f\in\mathcal{C}^\infty(\SSS^1)$ supported away from 1 and $U$ is a unitary operator on a Hilbert space $\mathcal{H}.$\\
$f\circ\psi$ is a smooth function on $\R$ with compact support. By Proposition \ref{complex} there exists an almost analytic 
extension $(f\circ\psi)^{\C}$ of $f\circ\psi,$ such that $(f\circ\psi)^{\C}$ supported in $\{x+i y;x\in supp (f\circ\psi), |y|\leq C<1\}.$\\

  Let $f^{\C}_{\SSS^1}=(f\circ\psi)^{\C}\circ\psi^{-1}$ defined on $\C.$ It is a smooth function, supported in\\
  $\Omega_0:=\psi(supp(f\circ\psi)\times[-C,C])$.\\
For $\xi\in \SSS^1, f^{\C}_{\SSS^1}(\xi)=(f\circ\psi)^{\C}\circ \psi^{-1}(\xi)=f\circ\psi\circ \psi^{-1}(\xi)=f(\xi),$ since $\psi^{-1}(\xi)\in\R.$\\
As in the proof of Theorem \ref{Sjos} we shall show that  (\ref{int}) converges.\\
As $ f^{\C}_{\SSS^1}\in\mathcal{C}^{\infty}(\C)$ with compact support  in $\C$  then, for all $\xi\in\SSS^1,$

\begin{equation}
 f(\xi)=(2i\pi)^{-1}\int_{\C}\frac{\partial_{\bar{z}} f^{\C}_{\SSS^1}(z)}{z-\xi}\,dz\wedge d{\bar{z}}
\end{equation}
by (\ref{compactif}) with $\omega=\C$.\\
Making a call to  lemma  \ref{derive} with $h=\psi^{-1}$ and $g=(f\circ\psi)^{\C}$ and  to proposition \ref{complex} with $l=1,$ 
and to lemma \ref{majpsi}, we get the  following uniform convergence on $\SSS^1$,
$$\lim_{n\rightarrow+\infty}\bigg\|\int_{d(z,\SSS^1)>1/n}\partial_{\bar{z}}f^{\C}_{\SSS^1}(z)(z-\cdot)^{-1}\,dz\wedge d{\bar{z}}-f(\cdot)\bigg\|_{\infty,\SSS^1}=0.$$
 Therefore, by Proposition \ref{conv},
$$\int_{\C}\partial_{\bar{z}}f^{\C}_{\SSS^1}(z)(z-U)^{-1}\,dz\wedge d{\bar{z}}$$
converges  in  operator  norm in $\mathcal{B}(\mathcal{H})$ and we have (\ref{uniHJ}).
\end{preuve}

\section{Notation}
$\R$ is the set of real numbers.\\
$\C$ is the complex plane.\\
$\SSS^1$ is the unit sphere in $\C$.\\
$D(0,1)$ is the open disk with center at the origin and of radius $1$ in $\C$.\\
$\omega$ is an open set in the complex plane $\C.$\\
$\partial\omega$ is the boundary of $\omega$ in $\C$.\\
$Im(z)$ is the imaginary of $z$. \\
$\C^+=\{z\in\C ; Im(z)>0\}.$\\
$\C^-=\{z\in\C ; Im(z)<0\}.$\\
$\C\setminus{\overline{D(0,1)}}$ is the complementary closed of the disk with center 0 and radius 1 in $\C$.\\
$\mathcal{C}^{\infty}(\R)$ is the set of infinitely differentiable functions on $\R.$\\
$\mathcal{C}^1(\omega)$  is the set of continuously differentiable functions on $\Omega$.\\
$\mathcal{H}$ is a Hilbert space.\\
$\mathcal{B}(\mathcal{H})$  is the set of bounded operators on $\mathcal{H}.$\\
$I$  is the identity operator on $\mathcal{H}.$\\
$\|A\|_{\mathcal{B}(\mathcal{H})}=sup_{\{u\in\mathcal{H};u\neq0\}}\frac{\|Au\|}{\|u\|}.$ \\
$\|f\|_{\infty,K}=sup_{x\in K}|f(x)|.$\\
$dz\wedge d{\bar{z}}$ is the Lebesgue measure on the complex plane.\\
$\partial_z=\frac{\partial_x-i \partial_y}{2}.$\\
$\partial_{\bar{z}}=\frac{\partial_x+i \partial_y}{2}.$\\
$\varphi^{(k)}$ denotes the kth derivative of $\varphi$.\\

\end{document}